\documentclass{amsart}
\usepackage{graphicx} 
\usepackage{booktabs}
\usepackage{multirow}
\usepackage{bm}
\usepackage[table]{xcolor}
\usepackage{braket}
\usepackage{todonotes}
\usepackage{subcaption}
\usepackage[foot]{amsaddr}
\usepackage{hyperref}

\title[Nonlinear Reduction strategies for Data Compression]{Nonlinear reduction strategies for data compression: a comprehensive comparison from diffusion to advection problems}

\author{Isabella Carla Gonnella$^1$, Federico Pichi$^{1}$, Gianluigi Rozza$^1$}
\address{$^1$ mathLab, Mathematics Area, SISSA, via Bonomea 265, I-34136 Trieste, Italy}

\date{}

\begin{document}

\begin{abstract}
  This work presents an overview of several nonlinear reduction strategies for data compression from various research fields, and a comparison of their performance when applied to problems characterized by diffusion and/or advection terms. We aim to create a common framework by unifying the notation referring to a common two-stage pipeline. At the same time, we underline their main differences and objectives by highlighting the diverse choices made for each stage. We test the considered approaches on three test cases belonging to the family of Advection-Diffusion problems, also focusing on the pure Advection and pure Diffusion benchmarks, studying their reducibility while varying the latent dimension. Finally, we interpret the numerical results under the lens of the discussed theoretical considerations, offering a comprehensive landscape for nonlinear reduction methods for general Advection-Diffusion dynamics.
\end{abstract}

\maketitle

\section{Introduction}

In recent years, the topic of \textit{Data Reduction} has become pervasive in the scientific computing community, as it contributes to popular ideas spread across many related fields, e.g.\ Numerical Analysis, Manifold Learning, Representation Learning, Generative Modeling, and more. To pursue a comprehensive analysis of the topic, we start by defining what is meant here by \textit{reduction}: identifying a low dimensional set of coordinates in the data space that encodes most of the information contained in the full dataset. The main assumption behind a reduction task is that the data manifold $\mathcal{M}$ at hand  lives, with a small approximation error, on a lower dimensional manifold $\mathcal{Z}_N$ of 
\begin{figure}[hbt!]
\centering
    \includegraphics[trim={0cm 20cm 0cm 2cm}, clip, width=1\textwidth]{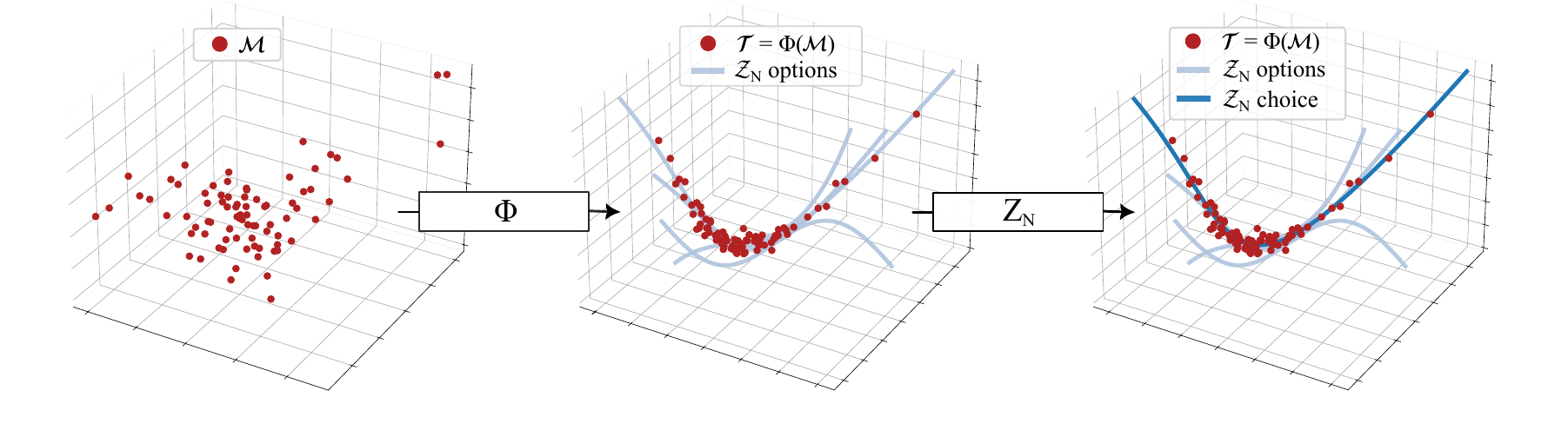}
    \caption{Nonlinear reduction two-stage pipeline: considering the original manifold $\mathcal{M}$, the map $\Phi$ transforms the data into the new manifold $\mathcal{T}$, on which reduction is performed through the operator $Z_N$.}
    \label{fig:general_reduction_scheme}
\end{figure}
dimension $N$ embedded in the higher dimensional space of dimension $D$, with $N\ll D$. Identifying such reduced space is convenient for interpretability but also helps speed up demanding computations, as it offers a suitable lower-dimensional dataset to process.

This work aims to uncover some important relations between reduction strategies developed in different field, trying to align the main perspectives by underlying their common ground and identifying their differences. Here, we want to address the characterization of nonlinear reduction strategies, generally adopted as a way to overcome the typical limitations of linear reduction when addressing highly nonlinear manifolds. To present a broad overview of different nonlinear approaches, we refer to the two-stage scheme reported in Figure \ref{fig:general_reduction_scheme} interpreting them as decomposable in two main stages, shaping their differences only by adopting different choices for the maps $\Phi$ and $Z_N$ and the reduced space $\mathcal{Z}_N$.

In particular, in Figure \ref{fig:general_reduction_scheme}, we represent the initial data-manifold $\mathcal{M}$ undergoing a series of transformations. The first one is given by the map $\Phi$ such that $\mathcal{T}=\Phi(\mathcal{M})$, and stands for a generic \textit{pre-processing step} that adjusts the original manifold to make it more suitable to the actual \textit{reduction step}, represented by the map $Z_N$ in the right part of the scheme. This map brings the transformed data manifold $\mathcal{T}$ into a low-dimensional space $\mathcal{Z}_N$ spanned by $N$ ``bases" that well capture the structure of $\mathcal{T}$, and associates each datum with its low-dimensional coefficients, for instance belonging to  $\mathbb{R}^N$, which characterize its reduced representation with respect to the chosen bases. The bright blue line precisely pictures the space $\mathcal{Z}_N$, embedded in the space of $\mathcal{T}$, that has been chosen between an ensemble of possible $N$-dimensional options as the one that best represents the (transformed) data.

\begin{figure}[b!]
\centering
    \includegraphics[trim={0cm 17cm 0cm 3cm}, clip, width=1\textwidth]{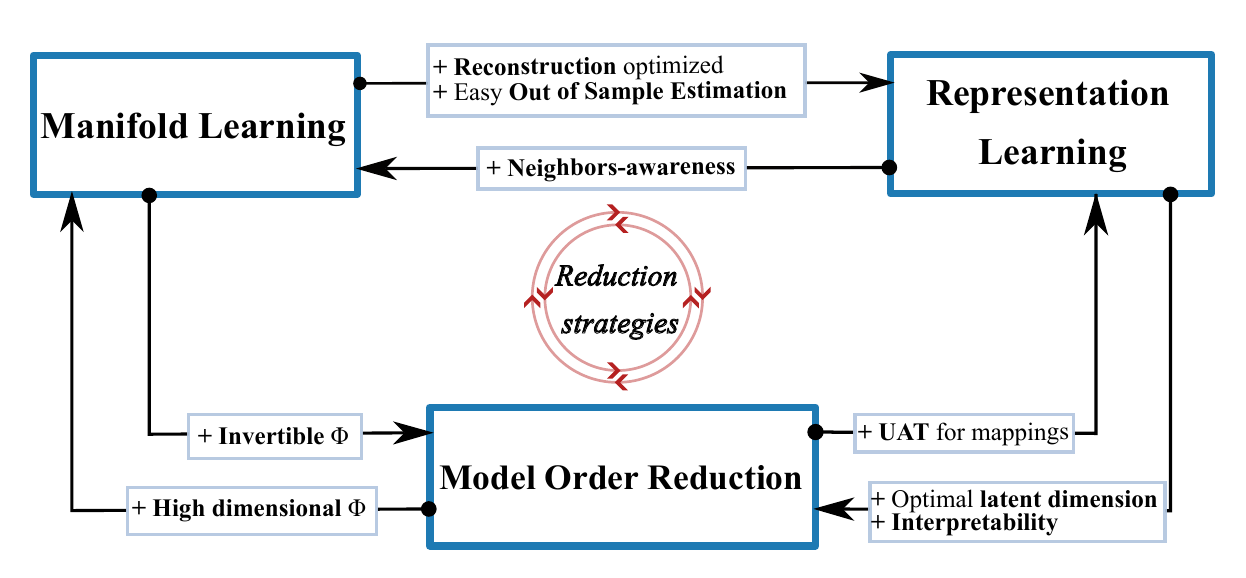}
    \caption{Bi-directional flow comparing the properties gathered respectively with methods of Manifold Learning, Representation Learning, and Model Order Reduction in the context of nonlinear strategies for data compression.}
    \label{fig:link_scheme}
\end{figure}

In this work, we address some well-known nonlinear reduction methodologies developed in the context of Numerical Analysis with Model Order Reduction techniques, Manifold Learning, and Representation Learning. Each one of them has some traditional strategies to perform data reduction, therefore, to give a glance at their connection and strengths, we report in Figure \ref{fig:link_scheme} a general scheme representing their links. Precisely, a double flow is proposed to report the properties gained in the reduction process passing from one methodology to another. These comments and connections will be deeply described and clarified in the following sections, nevertheless, it is worth giving a short brief overview of their main differences. Representation Learning usually builds the reduction problem as a supervised optimization of deep networks having strong approximation capabilities demonstrated in the Universal Approximation Theorem (UAT) \cite{kvalheim2024autoencoderswork}, aiming at minimizing the reconstruction error of each datum independently from its neighbors.  Manifold Learning learns a latent space by exploiting some data-structure awareness to properly cluster the data in the new low-dimensional latent space. Finally,  Reduced Order Models (ROMs) in Numerical Analysis, generally design an interpretable reduction pipeline providing convergence estimates on its reduced representations. On the contrary, learning-based algorithms usually lack interpretability, not clarifying the meaning of the obtained latent spaces, i.e.\ not providing information on the chosen bases for $\mathcal{Z}_N$ corresponding to the low-dimensional coefficients provided as latent representations of the data.

We start our analysis by introducing the main concepts and presenting \textit{linear reduction}, together with its main difficulties when dealing with some specific problems, for which we report a case test in Section \ref{sec:linear_reduction}. In particular, we refer to the setting of Advection-Diffusion problems, which notoriously model various physical phenomena such as heat and mass transfer, fluid motion, oceanography, air pollution, meteorology, and so on \cite{thongmoon2006comparison,jena2021computational}. Moreover, by parametrizing both the Advection and the Diffusion contributions in the differential model, one can obtain a set of significantly different dynamics to test reduction strategies on: \textit{pure Advection}, \textit{pure Diffusion}, and \textit{Advection-Diffusion} problems. Successively, in Section \ref{sec:different_same}, we dive into the details of the two-stage characterization for well-known ROM approaches (see Subsection \ref{subsec:numerical}), Manifold Learning (see Subsection \ref{subsec:manifold_learning}) and Representation Learning (see Subsection \ref{subsec:autoencoders}), pointing out in each case the construction strategies adopted for $\Phi, Z_N$ and $\mathcal{Z}_N$. Finally, in Section \ref{sec:results}, we report some numerical results regarding the application of the aforementioned reduction strategies to Advection-Diffusion problems, highlighting the pros and cons of each approach to the different simple, yet meaningful, test cases presented.
\section{Challenges in linear reduction}\label{sec:linear_reduction}
Given the parametric data manifold $\mathcal{M}:=\{u_{\mu_1},\dots,u_{\mu_n}\}$, where $u_{\mu_i}\in\mathcal{X}, \   \forall \mu_i\in\mathcal{P}:=\{\mu_1,\dots,\mu_n\}$ with $\mathcal{X}$ Hilbert space, traditional approaches to reduction rely on linear approximation spaces of the form $\mathcal{Z}_N=\text{span}\{\psi_j\}_{j=1}^N$ with $\psi_j\in\mathcal{X}, \ \forall j\in\{1,\dots, N\}$, so that

\begin{align*}
    u_{\mu_i} \simeq \Pi_{\mathcal{Z}_N}u_{\mu_i}:= \overline{u}+\sum_{j=1}^Nz_{j}(\mu_i)\psi_j,\ \forall \mu_i\in\mathcal{P},
\end{align*}

given the parameter-dependent reduced variables $z_j:\mathcal{P}\rightarrow\mathbb{R},\ \forall j\in\{1,\dots,N\}$ and $\overline{u}:=\frac{1}{n}\sum_{i=1}^nu_{\mu_i}$.

A well-known linear reduction method is the Principal Component Analysis (PCA) or Proper Orthogonal Decomposition (POD), depending on the field of application. It is particularly important because it offers the best-reduced representation of $\mathcal{M}$ in terms of $L_2$ error. Indeed, it uses as $\{\psi_j\}_{j=1}^N$ the $N$ orthogonal directions that maximize the variance of the data manifold $\mathcal{M}$. It can be proven that they are the first $N$ eigenfunctions of the covariance operator $\mathcal{K}:\mathcal{X}\rightarrow\mathcal{X}$ defined as
\begin{align*}
    \mathcal{K}(v):=\frac{1}{n}\sum_{i=1}^n\langle u_{\mu_i}-\overline{u}, v\rangle_{\mathcal{X}}(u_{\mu_i}-\overline{u})\quad\forall v\in\mathcal{X},
\end{align*}
so that
\begin{equation}\label{eq:pca_bases}
    \mathcal{K}(\psi_j)=\lambda_j\psi_j,\quad \forall j\in\{1,\dots,N\}.
\end{equation}
Consequently, the coefficients $\{z_j(\mu_i)\}_{j=1}^N$ are obtained by projecting the data on the new basis:
\begin{align*}
    z_j(\mu_i):=\langle u_{\mu_i},\psi_j\rangle_{\mathcal{X}},\quad\forall i\in\{1,\dots,n\},\ \forall j\in\{1,\dots,N\}.
\end{align*}

For instance, in the case of $\mathcal{X}\equiv\mathbb{R}^D$, if one defines the matrix of centered data as $U\in\mathbb{R}^{n\times D}$ (i.e.\ the mean of each column is zero), then $\{\psi_j\}_{j=1}^N$ are exactly the first $N$ eigenvectors of the covariance matrix $K=U^TU$.

Linear reduction thus relies on the availability of low-dimensional accurate representation of the manifold $\mathcal{M}$ itself. Nevertheless, this is not always the case, as the linear reducibility of $\mathcal{M}$ can be investigated via the Kolmogorov $N$-width $d_N(\mathcal{M})$:
\begin{align*}
d_N(\mathcal{M}):=\inf_{\mathcal{Z}_N}\sup_{u\in\mathcal{M}}\|u-\Pi_{\mathcal{Z}_N}u\|_{\mathcal{X}},
\end{align*}
where $\Pi_{\mathcal{Z}_N}:\mathcal{M}\rightarrow\mathcal{Z}_N$ is the orthogonal projection operator onto the reduced space $\mathcal{Z}_N$.  In particular, the Kolmogorov $N$-width can be seen as the worst approximation error that can be made on data belonging to $\mathcal{M}$ using the best linear subspace $\mathcal{Z}_N$ of dimension $N$ in terms of $\|\cdot\|_{\mathcal{X}}$. Furthermore, a good approximation of the decay of $d_N{(\mathcal{M})}$ while increasing $N$ is notoriously provided by the decay of the eigenvalues of the covariance operator $\mathcal{K}$ used in PCA \cite{floater2021best, arbes2024kolmogorovnwidthlineartransport,GreifDecayKolmogorovNwidth2019}.

While under certain circumstances, depending on the properties of $\mathcal{M}$, it is possible to demonstrate that $d_N(\mathcal{M})$ rapidly decays for increasing $N$ \cite{bachmayr2017kolmogorov, buffa2012priori}, there exists a large class of problems for which linear approximation does not offer an accurate reduced representation for a reasonably small amount ``modes". It is, for instance, the case of data collected from Advection-dominated dynamics \cite{TaddeiReducedBasisTechniques2015,NoninoReducedBasisMethod2023}, which are commonly described as the prototypical linearly nonreducible $\mathcal{M}$, i.e.\ characterized by a slow decay of $d_N(\mathcal{M})$.

\begin{figure}[hbt!]
    \subcaptionsetup[figure]{skip=-40pt,slc=off,margin={0pt,0pt}}

    \centering\subcaptionbox{\label{subfig:snapshots}}{\includegraphics[trim={0cm 0cm 0cm 0cm}, clip, width=0.47\textwidth]{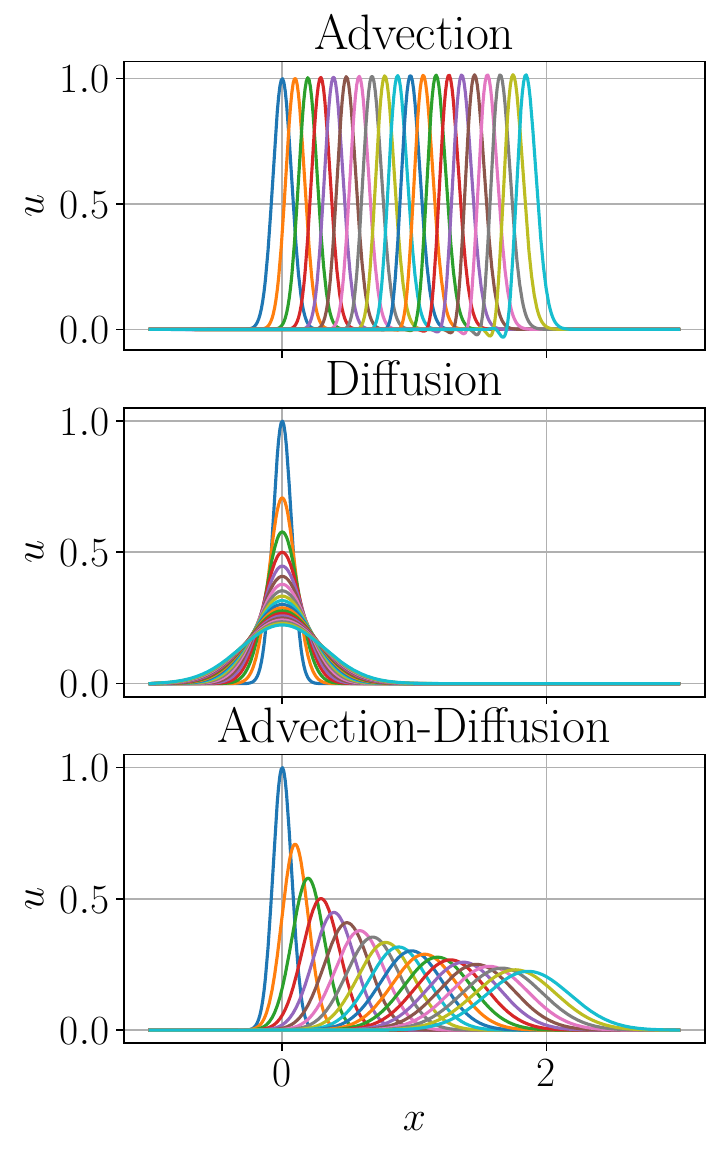}
    }
    \hfill
    \centering\subcaptionbox{\label{subfig:pca_eigen}}{\includegraphics[trim={0cm 0cm 0cm 0cm}, clip, width=0.47\textwidth]{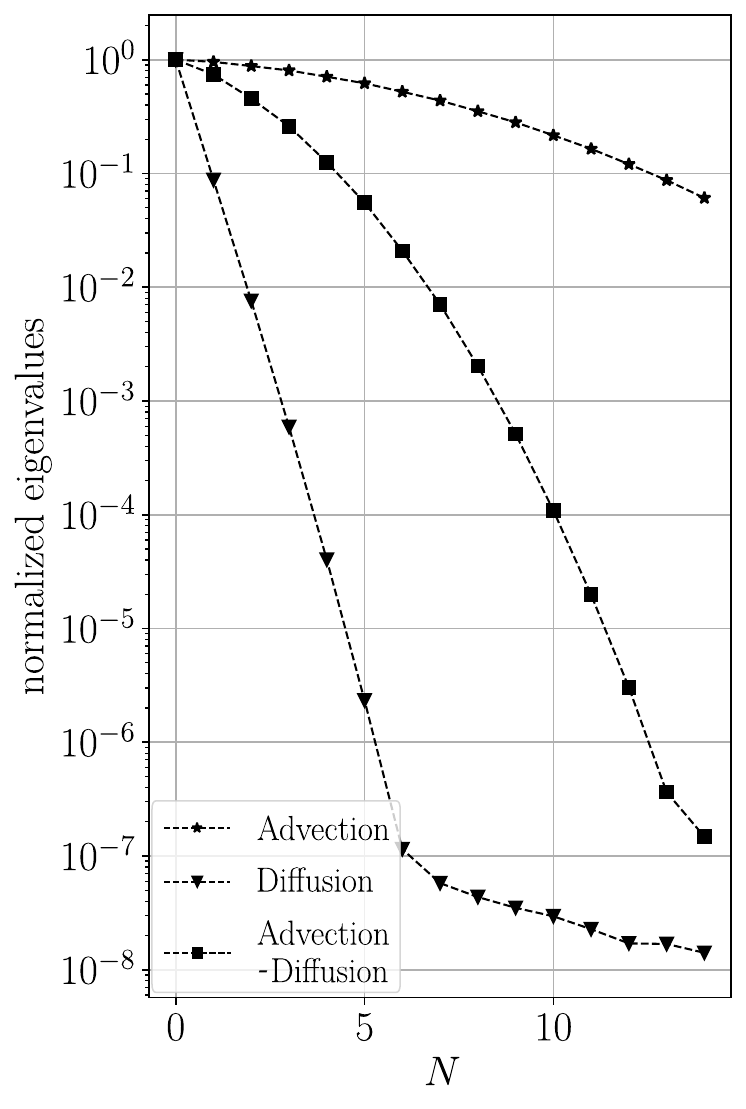}
    }

    \caption{On the left, some solutions at different times of the pure Advection ($c_T=4, c_D=0$), pure Diffusion ($c_T=0, c_D=0.1$), and Advection-Diffusion systems ($c_T=4, c_D=0.1$). On the right, their eigenvalues decay for increasing values of $N\in\{1,\dots,15\}$.}
    \label{fig:adv_diff_pca}
\end{figure}

This can be seen by comparing the decay of PCA eigenvalues of the manifold $\mathcal{M}$ obtained by solving the following Advection-Diffusion system for different values of the parameters $c_T$ and $c_D$ in the following equation in the spatial domain $\Omega\subset\mathbb{R}^D$:
\begin{equation}\label{eq:adv_diff}
    \begin{cases}
        \partial_t u(\bm{x},t) +c_T\nabla_{\bm{x}}\cdot u(\bm{x},t)-c_D\Delta_{\bm{x}} u(\bm{x},t)=0, &\quad\bm{x}\in\Omega,\ t>0\\
        u(\bm{x},0) = u_0,&\quad \bm{x}\in\Omega.
    \end{cases}
\end{equation}
Indeed, for $c_T=0, c_D\neq 0$ one obtains a purely diffusive system, while for $c_T\neq 0, c_D=0$ a purely advective one. 
Note that, when dealing with time-dependent systems, one may choose to include the time variable in the parameters' set. In Figure \ref{fig:adv_diff_pca} we compare the eigenvalues decay of the Advection, Diffusion, and Advection-Diffusion manifolds, both obtained by collecting solutions of Equation \eqref{eq:adv_diff} for $t\in \mathcal{P} = [0, T]$ with $T=0.5$ for fixed values of the other parameters: respectively $(c_T=4, c_D=0)$, $(c_T=0, c_D=0.1)$, and $(c_T=4, c_D=0.1)$. Thus, it is clear that, even in the ``non-parametric setting'' (only evolving the systems in time), linear reduction is not well-suited for Advection-dominated problems, since the correspondent PCA eigenvalues exhibit a significantly slower decay w.r.t.\ the other models, for which diffusive effects are present.

\section{Nonlinear Reduced Models: different Learning Strategies, same Reduction Pipeline} \label{sec:different_same}

In this section, we address well-known nonlinear reduction strategies from different research fields by linking them to the scheme reported in Figure \ref{fig:general_reduction_scheme}. In particular, we identify a nonlinear reduction process with a two-stage pipeline composed of a possible \textit{pre-processing} step transforming the manifold $\mathcal{M}$ into $\mathcal{T}$ through the operator $\Phi:\mathcal{M}\rightarrow\mathcal{T}$, and of a standard \textit{reduction step}, which aims to find the operator $Z_N:\mathcal{T}\rightarrow\mathcal{Z}_N\subseteq\mathbb{R}^N$ such that the 
reduced manifold $\mathcal{Z}_N:= Z_N(\mathcal{T})$ has dimension $N$ and approximates sufficiently well $\mathcal{T}$. Thus, the main goal of the first step is to obtain a ``simpler" representation of the manifold $\mathcal{T}$ which is more suitable for reduction, while the second one actually searches for the more convenient $N$-dimensional description of the resulting manifold $\mathcal{T}:=\Phi(\mathcal{M})$. 
For instance, in the case of projection-based ROMs, the last step aims to find $N$ bases $\{\psi_j\}_{j=1}^N$ on which to project $\mathcal{T}$, so that $Z_N$ is taken as a projection operator. 
In addition, since it is generally desirable to reconstruct the original data $\mathcal{M}$ from $\mathcal{Z}_N$, the methods should also provide the operators $\Phi^{-1}:\mathcal{T}\rightarrow\mathcal{M}$ and $Z_N^{-1}:\mathcal{Z}_N\rightarrow\mathcal{T}$. Nevertheless, we underline the abuse of notation for $\Phi^{-1}$ and $Z_N^{-1}$, as they are defined only to provide a way to map back the reduced representation to the original manifold, being $\Phi$ and $Z_N$ not necessarily invertible operators. Therefore, each method gives its own characterization to $\Phi^{-1}$ and $Z_N^{-1}$ (see for instance Section \ref{subsec:autoencoders} for an example of non-invertible $Z_N$).

Recalling Section \ref{sec:linear_reduction}, in the linear PCA/POD approach one can identify $\Phi$ with the identity operator, so that $\mathcal{T}=\mathcal{M}$, and $Z_N:\Phi(u_{\mu_i})\mapsto \{z_j(\mu_i)\}_{j=1}^N\in\mathcal{Z}_N$ with the operator assigning to each transformed datum the reduced variables $\bm{z}(\mu_i):=\{z_j(\mu_i)\}_{j=1}^N$ for a given parameter value $\mu_i\in\mathcal{P}$ corresponding to the projection of $\Phi(u_{\mu_i})$ on the bases $\psi_j$. Then the ``inverse'' operator $Z_N^{-1}(\bm{z})=\sum_{i=j}^Nz_j\psi_j$ projects back the reduced variables on $\mathcal{T}$, being $\{\psi_j\}_{j=1}^N$ defined as in \eqref{eq:pca_bases} and corresponding to the chosen bases of $\mathcal{T}(\equiv\mathcal{M})$.
Moreover, we highlight that the same reasoning applies for methods that rely on the interpolation of different local projection-based ROMs correspondent to different parameters values \cite{amsallem2011online, peherstorfer2014localized}. Indeed, such approaches can be represented by taking $\Phi$ as the identity map (similarly to the classical PCA case), and choosing $Z_N$ to be a projection-operator on a set of parameter dependent bases $\{\psi_j(\mu)\}_{j=1}^N$, so that $\mathcal{Z}_N\equiv\mathcal{Z}_N(\mu)$.

Now, we continue by discussing how different nonlinear reduction techniques fit in the above pipeline. With this aim, we present below the characterization of these two stages for some numerical nonlinear ROMs, for well-known Manifold Learning techniques, and for Auto-Encoders, belonging to the Representation Learning field.

\subsection{Techniques from Model Order Reduction: Registration methods}\label{subsec:numerical}

Nonlinear approaches to reduction, as underlined in \cite{taddei2020registration,torlo2020model}, can be subdivided into \textit{Eulerian} and \textit{Lagrangian} ones. Eulerian methods search for possibly nonlinear parameter-dependent maps $Z_N\equiv Z_{N,\mu}$ and $Z_N^{-1}\equiv Z_{N,\mu}^{-1}$, with $\mu\in\mathcal{P}$, while keeping $\Phi\equiv Id$. Lagrangian approaches instead aim to perform classical linear reduction on $\mathcal{T}=\Phi(\mathcal{M})$ with a general nonlinear map $\Phi$, obtained by properly optimizing the image set $\mathcal{T}$ to make it well-suited for linear reduction \cite{NoninoReducedBasisMethod2023}.

Concerning Eulerian approaches, examples might rely on Grassmannian learning \cite{zimmermann2018geometric}, displacement interpolation \cite{rim2018model,CucchiaraModelOrderReduction2023}, or transformed snapshot methods \cite{nair2019transported, reiss2018shifted,RimManifoldApproximationsTransported2020}. Nevertheless, in this section, we are focusing on Lagrangian methods, given their capability to exploit the reduction pipeline of Figure \ref{fig:general_reduction_scheme}.

Regarding Lagrangian nonlinear approaches, we start by reporting the example of \textit{Registration Methods}, which generalizes the linear Reduced Basis (RB) technique \cite{hesthavenCertifiedReducedBasis2015,QuarteroniReducedBasisMethods2016,RozzaRealTimeReduced2024} by coupling it with a previous transformation of the data space. 
The term \textit{registration} has been used in the fields of computer vision and pattern recognition, to refer to the process of finding a spatial transformation that aligns two datasets.  In the context of the Numerical Analysis of parametric Partial Differential Equations (pPDEs) it has been adopted to indicate  the process of finding a transformation $\Phi$ that improves the linear reduction of a given manifold exploiting such alignment strategy.


Indeed, \textit{Registration Methods} \cite{ferrero2022registration, taddei2021space, taddei2022optimization} deal with the task of finding a parameter-dependent bijection $\Phi\equiv\Phi_{\mu}$ such that, given the initial manifold $\mathcal{M}$ of sampled solutions of a pPDE $\{u_{\mu_i}(\bm{x})\}_{i=1}^n$ with $\bm{x}\in\Omega\subset\mathbb{R}^D$, for an ensemble of different values of the parameter $\{\mu_i\}_{i=1}^n$, the image manifold $\mathcal{T}:=\Phi(\mathcal{M})$ is suitable for linear reduction, i.e.\ has a nice Kolmogorov decay. This is done by defining the \textit{optimal map} $\Phi_{\mu}$ as the one that minimizes the $L_2$ distance between the transformed data and the reference solution $u_{\overline{\mu}}\in\mathcal{M}$, for which $\Phi_{\overline{\mu}}$ is the identity map. In particular, $\Phi_{\mu}$ is searched as an affine map of the form
\begin{align*}
    \Phi_{\mu}^{\bm{\hat{a}}}(\bm{x}) := \bm{x}+\sum_{m=1}^M \hat{a}_m(\mu)\varphi_m(\bm{x}),
\end{align*}
where $\{\varphi_m\}_{m=1}^M$ is a given set of bases that well approximate the space $L_2(\Omega)$ (e.g.\ the set of Legendre polynomials), and $\hat{\bm{a}}=\{\hat{a}_m\}_{m=1}^M$ are coefficients obtained through the following regularized constrained minimization problem, where $\varepsilon, C$ serve to weakly impose $J(\Phi_{\mu}^{\bm{\hat{a}}}(\bm{x}))\in[\varepsilon,\frac{1}{\varepsilon}]$ for the Jacobian, and thus that $\Phi_{\mu}^{\bm{\hat{a}}}(x)$ is a bijection from $\Omega$ into itself:

    \begin{equation}
    \begin{split}\label{eq:registration_optimization}
     \hat{\bm{a}}=\arg&\min_{\bm{a}\in\mathbb{R}^M} \|u_{\mu}\circ\Phi^{\bm{a}}_{\mu}(\bm{x})-u_{\overline{\mu}}(\bm{x})\|_{L_2(\Omega)}^2 + \xi|\Phi_{\mu}^{\bm{a}}(\bm{x})|^2_{H_2(\Omega)} \\&\text{s.t. }\int_{\Omega}e^{\left(\frac{\varepsilon-J(\Phi_{\mu}^{\bm{\hat{a}}}(\bm{x}))}{C}\right)}+e^{\left(\frac{J(\Phi_{\mu}^{\bm{\hat{a}}}(\bm{x}))-\frac{1}{\varepsilon}}{C}\right)}d\bm{x}\leq\delta. 
     \end{split}
     \end{equation}
The $H_2$ regularizing term measures the deviations from linear maps, as for $v\in H_2(\Omega)$ linear the seminorm $|v|_{H_2(\Omega)}:=\sum_{i,j=1}^D\int_{\Omega}(\partial^2_{i,j}v)^2=0$. This can be further interpreted as a Tikhonov regularization and has the effect of controlling the gradient of the Jacobian, which is important for ensuring invertibility of $\Phi_{\mu}^{\bm{\hat{a}}}(x)$.

Thus, Registration Methods offer a general and interpretable framework, clarifying the optimization process that leads to a given choice of $\Phi$ while ensuring its invertibility and the possibility of reconstructing the original samples. Furthermore, since they rely on classical linear methods for the reduction step, they can be coupled with any reduction strategy by selecting the dimension $N$ looking at the eigenvalues decay of the covariance operator computed on the transformed manifold $\mathcal{T}=\Phi(\mathcal{M}):=\{u_{\mu}\circ\Phi_{\mu}^{\hat{a}}\ |\ u_{\mu}\in\mathcal{M}\}$, not imposing it a-priori.

Moreover, an interesting direction for Registration has been offered by Optimal Transport (OT) theory, which focuses on the existence and computability of \textit{optimal transport maps} between two probability measures with finite second order moments $\rho,\sigma\in\mathcal{P}_2(\Omega,\mathcal{W}_2)$ with respect to the Wasserstein metric $\mathcal{W}_2$. It can be reinterpreted from a Registration perspective by considering $u_{\mu_i}(\bm{x})\in\mathcal{P}_2(\Omega,\mathcal{W}_2), \ \forall\mu_i\in\mathcal{P}$, computing $\Phi$ through the transport maps $T_{u_{\mu}\rightarrow u_{\overline{\mu}}}$ between the distributions $u_{\mu}(\bm{x})$ and the reference element $u_{\overline{\mu}}(\bm{x})$. Some works developed and discussed this intuition, proposing frameworks in which OT maps can be useful to find a proper $\mathcal{T}=\Phi(\mathcal{M}):=\{(T_{u_{\mu}\rightarrow u_{\overline{\mu}}})_{\#}u_{\mu}\ |\ u_{\mu}\in\mathcal{M}\}$ in case of problems with slow decay, such as the Advection-dominated ones \cite{blickhan2023registration, ehrlacher2020nonlinear}.

In particular, \cite{blickhan2023registration} proposes a way to build a proper transformation from the manifold $\mathcal{M}\subset\mathcal{P}(\Omega)$ to a space as simple as possible, even coinciding with a proper reference solution $u_{\overline{\mu}}$ for a given reference value of the parameter $\overline{\mu}$. In this setting, $\Phi_{\mu}$ coincides with the OT map $T_{u_{\mu}\rightarrow\overline{u}}$ between $u_{\mu}$ and $u_{\overline{\mu}}$, as it minimizes the following problem
\begin{align*}
    T_{u_{\mu}\rightarrow u_{\overline{\mu}}} = \arg\min_{T}\int_{\Omega}\frac{1}{2}|T(u_{\mu})-u_{\overline{\mu}}|^2d\bm{x},
\end{align*}
thus corresponding to the \textit{optimal map} according to the definition of Wasserstein distance, which reads
\begin{align*}
    \mathcal{W}_2(u_{\mu},u_{\overline{\mu}}) = \min_{T}\int_{\Omega}\frac{1}{2}|T(u_{\mu})-u_{\overline{\mu}}|^2d\bm{x}.
\end{align*}

\subsection{Techniques from Manifold Learning: Kernel and Spectral methods} \label{subsec:manifold_learning}
Manifold Leaning belongs to the Metric Learning field, which, by definition, aims at optimizing a problem-dependent metric that provides a well-suited similarity measure to represent the given dataset by means of a few \textit{clusters}. This task can be approached both in the supervised and unsupervised context, and in the latter case it precisely takes the name of Manifold Learning. Such techniques aim to understand the manifold geometrical and topological features. Precisely, this is often done by identifying the bases spanning an informative lower dimensional space, on which to project the original data manifold to both preserve and emphasize possibly hidden core information. As a result, manifold learning often provides significant hints for the reduction, as it searches for a low dimensional representation of the high dimensional data while preserving relevant information of its original topology \cite{meilua2024manifold, ma2012manifold, kaya2019deep,bellet2015metric}.

For instance, PCA (see Section \ref{sec:linear_reduction}) belongs to this class of methods, as it finds new bases maximizing the variance of the reduced representation so that the original data can be treated with a reduced number of coordinates while still preserving its intrinsic variance. A nonlinear generalization to PCA is offered by its ``kernelized" version \cite{scholkopf1998nonlinear, blanchard2007statistical}, which only modifies the $\Phi$. Indeed, in kernel PCA (kPCA), $\Phi$ is assumed to be a high-dimensional function defining the correspondent kernel $k$ used in computations. The idea behind a high-dimensional $\mathcal{T}$ is to obtain a sparsification of data on new directions advantageous to the subsequent PCA reduction.

To enable the computation of the PCA decomposition of the high-order nonlinear transformation of the original manifold, the so-called \textit{kernel trick} is introduced. More precisely, if each datum $u_{\mu_i}$ belongs to $\mathcal{X}\equiv\mathbb{R}^D$ so that $\bm{U}:=\{u_{\mu_i}\}_{i=1}^n\in\mathbb{R}^{n\times D}$ and $\Phi:\mathbb{R}^D\rightarrow\mathbb{R}^M$ with $M$ arbitrarily large, the kernel is thus defined as a dot product in the \textit{feature space}  $\mathbb{R}^N$
\begin{align*}    k(u_{\mu_i},u_{\mu_j})=\braket{\Phi(u_{\mu_i}),\Phi(u_{\mu_j})},
\end{align*}
and the \textit{kernel trick} enables to compute the kPCA reduction algorithm without explicitly knowing $\Phi$, but only relying on $k$.

Indeed, performing PCA on the \textit{feature space}, we study the spectral decomposition of the matrix $\bm{C}:=\bm{\Phi}^T\bm{\Phi}$, where $\bm{\Phi}\in\mathbb{R}^{(n\times M)}$ is the evaluation matrix of the function $\Phi$ at the elements of the data manifold $\mathcal{M}=\{u_{\mu_i}\}_{i=1}^n$. It can be verified that, being $\{\lambda^{\bm{C}}_j\}_{j=1}^R$ and $\{\bm{v}^{\bm{C}}_j\}_{j=1}^R$ respectively the eigenvalues and eigenvectors of the matrix $\bm{C}$ of rank $R$
 \begin{align*}
     \bm{C}\bm{v}^{\bm{C}}_j = \lambda^{\bm{C}}_j \bm{v}^{\bm{C}}_j, \quad \forall j\in\{1,\dots,R\},
 \end{align*}
 they can be put in relation with the ones $\{\lambda^{\bm{K}}_j\}_{j=1}^R$ and $\{\bm{v}^{\bm{K}}_j\}_{j=1}^R$ of the \textit{kernel matrix} $\bm{K}_{i,j}:=k(u_{\mu_i},u_{\mu_j})=\bm{\Phi}\bm{\Phi}^T\in\mathbb{R}^{(n\times n)}$ which by definition has the same rank $R$. Thus, one obtains
 \begin{align*}
     \bm{v}_j^{\bm{C}} = \bm{\Phi}^T\bm{v}_j^{\bm{K}} \quad \text{and}\quad \lambda_j^{\bm{C}}=\lambda_j^{\bm{K}}\quad\forall j\in\{1,\dots,R\},
 \end{align*}
and the normalizing conditions $\lambda_j^{\bm{K}}(\bm{v}_j^{\bm{K}}\cdot\bm{v}_j^{\bm{K}})=1,\ \forall j\in\{1,\dots, R\}$, are applied to obtain $(\bm{v}_k^{\bm{C}}\cdot\bm{v}_k^{\bm{C}})=1$.

This way, the reduced representation of a given data $u\in\mathbb{R}^D$ can be identified only using the kernel $k$ without explicitly knowing $\Phi(u)$, but only computing its projection of the transformed onto the first $N<R$ principal components $\{\bm{v}_j^{\bm{C}}\}_{j=1}^N$:
 \begin{align*}
     \Phi(\bm{u})\bm{v}^{\bm{C}}_j =  \Phi(\bm{u})\bm{\Phi}^Tv_j^{\bm{K}}=\sum_{l=1}^n k(u,u_l)(\bm{v}^{\bm{K}}_j)_l,\quad\forall j\in\{1,\dots,N\}.
 \end{align*}

 Nevertheless, it is to be underlined that without $\Phi$ we can not reconstruct the original samples $u\in\mathcal{M}$ starting from their reduced representation. In particular, given a kernel $k$ satisfying some general assumptions, Mercer's theorem \cite{steinwart2012mercer} ensures the existence and uniqueness of $\Phi$, but it is in general too complicated to be explicitly written.

However, for the identification of relevant clusters of significantly nonlinear data manifolds, many different kernels have proved to be beneficial, and in this section we report some well-known Manifold Learning techniques that can be reconducted to kernelized versions of PCA \cite{ghojogh2023elements}, as it is the case of the ones listed in Table \ref{tab:kernel_methods}.

 \newcommand{\BB}[1]{\cellcolor{blue!#1}}
\renewcommand{\arraystretch}{1.2}
\begin{table}[htb!]
\centering
\begin{tabular}{l|c}
\rowcolor[HTML]{C1CDCD}
\toprule
\textbf{Method} & \textbf{Kernel} \\
\midrule
\textit{{Principal Component Analysis}} &\BB{20} $\bm{U}^T\bm{U}$  \\
\textit{{Multi-Dimensional Scaling}} &\BB{10}$-\frac{1}{2}\bm{H}\bm{D}\bm{H}$ \\
\textit{{Isomap}} &\BB{20} $-\frac{1}{2}\bm{H}\bm{D}^{(g)}\bm{H}$ \\
\textit{{Spectral Clustering}} & \BB{10} $\bm{L}^{\dag}$ \\
\textit{{Local Linear Embedding}} &\BB{20} $\bm{M}^{\dag}$ \\
\bottomrule
\end{tabular}
\vspace{.3cm}
\caption{Manifold Learning methods written as kPCA ones with different kernel functions.}
\label{tab:kernel_methods}
\end{table}

 \subsubsection{Multi-Dimensional Scaling (MDS)}

MDS 
\cite{ghojogh2023elements,saeed2018survey} aims to preserve a notion of distance between points, which is encoded in the matrix $\bm{D}$ and then \textit{double-centered} via the matrix $\bm{H}:=\bm{I}_n-\frac{1}{n}\bm{1}_n$. This is done by projecting the data on the directions which maximise the distance between the \textit{training points}, similarly to what is done with the matrix $K$ in the classical case of PCA, to preserve the variance in the reduced representation.

\subsubsection{Isomap}
It follows the same principles as MDS, only changing the distance used to compute matrix $\bm{D}^{(g)}$, which is not the Euclidean distance, but the geodesic one, defined as the length of the shortest path between two points on the
 possibly curvy (i.e., nonlinear) manifold. It can be computed by approximating it with a sum of local Euclidean distances between neighbour points \cite{choi2007robust,tenenbaum2000global,balasubramanian2002isomap}.

 \subsubsection{Spectral Clustering} The main idea behind Spectral Clustering (SC) is to divide the data $\bm{U}$ into clusters which are as separated as possible \cite{von2007tutorial,ng2001spectral}. It is done through an optimization procedure based on the so-called \textit{adjacency matrix} $\bm{W}$, defined as
 \begin{equation}\label{eq:adjacency}
     \bm{W}_{i,j} = \begin{cases}
         w(\bm{u}_i,\bm{u}_j)\quad &\bm{u}_j\in\text{Neighbours}(\bm{u_i}), \\
         0 \quad &\bm{u}_j\notin\text{Neighbours}, \\
         0 \quad &\bm{u}_i=\bm{u}_j,
     \end{cases}
 \end{equation}
accordingly to a \textit{weighting function} $w:\mathbb{R}^D\times\mathbb{R}^D\rightarrow\mathbb{R}$. The optimization process can be translated in a minimization of the sum of the weights of the points belonging to different clusters
\begin{align*}
    \min_{\bm{v}}\sum_{i=1}^n\sum_{j=1}^n w(\bm{u}_i,\bm{u}_j)(\bm{v}_i-\bm{v}_j),
\end{align*}
where $\bm{v}$ stands for a discriminator function for a specific cluster.  In this way, the reduced representation can be more suitable to be analyzed with simple clustering algorithms.

It has been proven that such problem is equivalent to the eigenvalue problem of the Laplacian matrix
\begin{align*}
    \bm{L}:=\bm{W}_{D}-\bm{W},
\end{align*} where $\bm{W}_D:=diag(\{d_1,\dots,d_n\})$ with $d_i:=\sum_{j=1}^nw_{i,j}$. Thus, this directly translates in the kPCA maximization problem taking the pseudo-inverse of the Laplacian $\bm{L}^{\dag}$ as \textit{kernel matrix}. Therefore, SC is equivalent to a kPCA version in which particular attention is given to the adjacency links of \textit{training data}, wanting to maximise the information encoded by the Laplacian. Giving an intuitive interpretation of $\bm{L}$, we remark that it can be seen as a measure of how \textit{curved} the neighborhood of a point is, so the SC approach aims at preserving the geometrical ``second order" characteristics of the points neighborhoods.

\subsubsection{Locally Linear Embedding}

Among the techniques searching for a reduced representation of $\bm{U}$, the goal of Locally Linear Embedding (LLE) methods is to preserve the linear reconstruction weights of data neighbors \cite{saul2000introduction,ghojogh2020locally, zhang2006mlle}. Indeed such approach is divided in two fundamental optimization steps: the first in the high-dimensional space, and the second in the low-dimensional one. During the first stage, one searches for the weight matrix $\bm{\tilde{W}}$ defined as
\begin{equation}\label{eq:LLE}
    \bm{\tilde{W}}:=\arg\min_{\bm{\tilde{W}}}\sum_{i=1}^n\left\|\bm{u}_i-\sum_{j\in\text{kNN}(\bm{x}_i)}\tilde{w}_{i,j}\bm{u}_j\right\|,
\end{equation}
so that a local linear reconstruction of each point is obtained. Then, the next stage exploits $\bm{\tilde{W}}$ to search for the optimal points $\bm{Y}$ in the reduced space, thus solving an second minimization problem of the form
\begin{align*}
    \bm{Y}:=\arg\min_{\bm{Y}}\sum_{i=1}^n\left\|\bm{y}_i-\sum_{j\in\text{kNN}(\bm{u}_i)}\tilde{w}_{i,j}\bm{y}_j\right\|,
\end{align*}
that can actually be seen as an eigenvalue problem for the matrix $\bm{M}$ defined as
\begin{align*}
    \bm{M}:=(\bm{I}-\bm{\tilde{W}})^T(\bm{I}-\bm{\tilde{W}}).
\end{align*}
Also in this case, in order to cast the problem as a kPCA one, the minimization can be transformed into a maximization problem by taking the pseudo-inverse $\bm{M}^{\dag}$.



\subsection{Techniques from Representation Learning: Auto-Encoders} \label{subsec:autoencoders}
From the field of Representation Learning, we include in the analysis Auto-Encoders (AE), one of the main Deep Learning (DL) architectures used for reduction \cite{MorrisonGFNGraphFeedforward2024a,PichiGraphConvolutionalAutoencoder2024,FrescaComprehensiveDeepLearningBased2021,FrancoLatentDimensionDeep2023,LeeModelReductionDynamical2020}. Indeed, Data Reduction has become a popular pre-processing step, prior to the actual learning algorithms, as it has been understood that the effectiveness and interpretability of machine learning methods largely relies on the selection of data features they are applied to. As a result, a significant part of the research work has been focusing on designing pipelines for data transformations that produce a suitable \textit{representation} for successful learning. To this aim, Auto-Encoders have been presented as a peculiar deep learning architecture to reduce the inputs dimension, in order to build an informative surrogate representation of them. In particular, Auto-Encoders are composed by an encoder and a decoder, which are layers that respectively reduce and reconstruct the original data manifold.

We can refer to the pipeline reported in Figure \ref{fig:general_reduction_scheme} also for their description. Indeed, they search for the reduced representation which \textit{minimize} the reconstruction error on the training dataset, relying on the UAT, that enables a representation through a deep network of a large variety of invertible encoding maps satisfying such minimization constraint. The fixed reduced dimension $N$ is expressed by construction by the number of neurons of the latent representation, and thus has to be established a priori to the training phase. In our notation, the encoder $f_e$ approximates the function $Z_N:\mathcal{X}\rightarrow\mathbb{R}^N$ that transforms the datum $u_{\mu_i}$ correspondent to the parameter $\mu_i$ in its latent representation $\bm{z}_{\mu_i}\in\mathbb{R}^N$, while the decoder $f_d$ corresponds to $Z_N^{-1}$ that maps back into $\mathcal{X}$ the latent representation of the input. 
Nevertheless, in contrast to the Registration and Manifold Learning cases, for Auto-Encoders we can not assume $\mathcal{Z}_N$ 
to be a linear space and $Z_N$ to be a projection operator with respect to some bases $\{\psi_j\}_{j=1}^N$, as it is nontrivial to explicitly identify such bases hypothetically spanning the latent space $Z_N$. Indeed, even assuming the encoder as a composition of a nonlinear map $\Phi$ and a reduction step $Z_N$, there is no easy way to characterize the latent space \cite{magri2022interpretability}.


Given the dataset $\bm{U}=\{\bm{u}_{\mu_i}\}_{i=1}^n$ with $\bm{u}_{\mu_i}\in\mathcal{X}\equiv\mathbb{R}^D,\ \forall i\in\{1,\dots,n\}$, the training procedure follows the minimization of the loss $\mathcal{L}$:
\begin{equation}
    \mathcal{L}(\bm{U}) := \frac{1}{n}\sum_{i=1}^n \|\bm{u}_i- f_d\circ f_e(\bm{u}_i)\|_{\mathcal{X}},
\end{equation}
where $f_e$ and $f_d$ are respectively denoting  the encoder and the decoder functions.

Additional terms can be inserted into the loss formulation as regularizing factors, as it is the case of Sparse Auto-Encoders \cite{ng2011sparse}, and of Contractive Auto-Encoders \cite{rifai2011higher}. The former class enforces sparsity in the hidden activations of the layers $\bm{A}$ inserting a regularized $L_1$ norm term:
\begin{equation}
    \mathcal{L}_{sparse}(\bm{U}) := \frac{1}{n}\sum_{i=1}^n \|\bm{u}_i- f_d\circ f_e(\bm{u}_i)\|_{\mathcal{X}} + \lambda\|\bm{A}\|_{L_1}.
\end{equation}
The second class instead aims to build an architecture resistant to input non-informative fluctuations, with the goal of inferring a stable latent representation. In order to do that, a regularizing term is inserted regarding the Jacobian matrix of the hidden layer $\bm{J}_A$, which represents the derivative of each hidden node with respect to each input component:
\begin{equation}
    \mathcal{L}_{contractive}(\bm{U}) := \frac{1}{n}\sum_{i=1}^n \|\bm{u}_i- f_d\circ f_e(\bm{u}_i)\|_{\mathcal{X}} + \lambda\|\bm{J}_A\|_{L_2},
\end{equation}
so that high derivatives w.r.t.\ inputs are discouraged.

Therefore, Auto-Encoders obtain a latent representation through loss minimizations that are not explicitly aware of the topological characteristics of the manifold, as instead happened in Section \ref{subsec:manifold_learning} with Manifold Learning techniques. On the other hand, such DL architecture focuses on enabling an accurate reconstruction of samples, which, in contrast, is not always possible in the previous cases due to the lack of an explicit $\Phi^{-1}$ operator. 

Furthermore, in many cases Auto-Encoders provide an advantageous reduction architecture for Out Of Sample (OOS) estimation, i.e.\ the identification of the reduced representation of data not belonging to the \textit{training set}. 
Indeed, for all the approaches presented in Sections \ref{subsec:numerical} and \ref{subsec:manifold_learning}, one has to perform some sort of multivariate regression strategy on the set of ``train reduced variables'' $\{\mu_i,\bm{z}_i\}_{i=1}^n$ to obtain the reduction operator $Z_N:u_{\mu}\mapsto\bm{z}_{\mu}\ \forall\mu\in\mathcal{P}$. In the DL case, instead, $Z_N$ is the encoder itself, which automatically maps any new input in a new point of the learnt reduced space.
\section{Numerical Results} \label{sec:results}

In this section, we apply the nonlinear reduction methods described above to the Advection-Diffusion problem described in Equation \eqref{eq:adv_diff} (see Section \ref{sec:linear_reduction}). In particular, we test them on the pure Advection problem ($c_D=0, c_T=4$), pure Diffusive problem ($c_D=0.1, c_T=0$), and Advection-Diffusion one ($c_D=0.1, c_T=4$). For simplicity, in each test case, the parameters are kept fixed and the time is evolved in the range $[0,T]$ with $T=0.5s$, and we fix $\Omega:=[-1,3]$. The \textit{training set} is composed by $20$ solution samples $u_{t_i}$ for $t_i$ uniformly distributed in $[0,T]$, while the \textit{testing set} by $200$ samples uniformly distributed in $[0,T]$ \footnote{The Python codes with the implementation of the methods exploited for the numerical results will be made available upon request.}.

\begin{figure}[t!]
    \subcaptionsetup[figure]{skip=-10pt,slc=off,margin={0pt,0pt}}
    \centering\subcaptionbox{\label{subfig:reg_eigenvalues_T}}{\includegraphics[trim={0cm 0cm 0cm 0cm}, clip, width=0.31\textwidth]{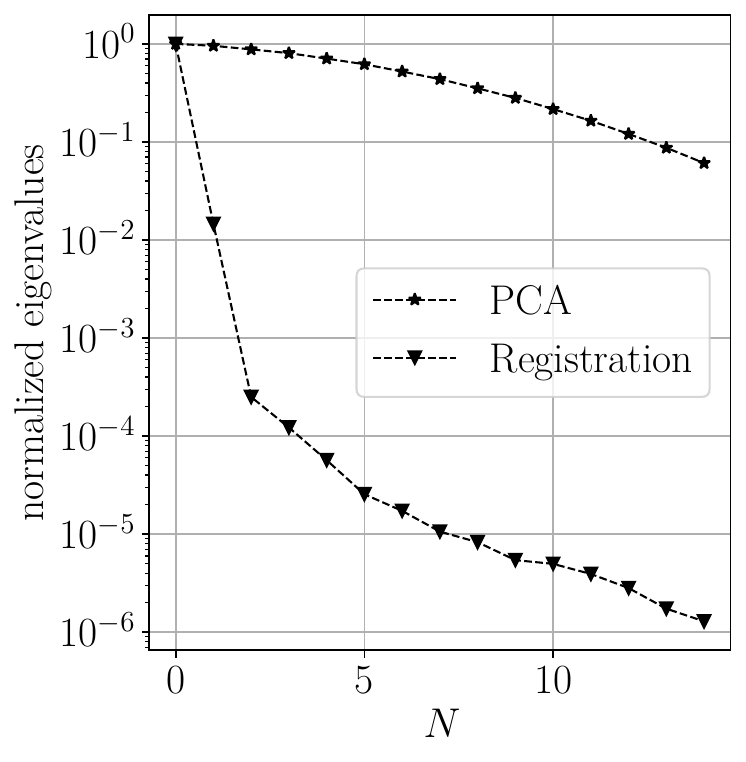}
    }
    \hfill
    \centering\subcaptionbox{\label{subfig:reg_eigenvalues_D}}{\includegraphics[trim={0cm 0cm 0cm 0cm}, clip, width=0.31\textwidth]{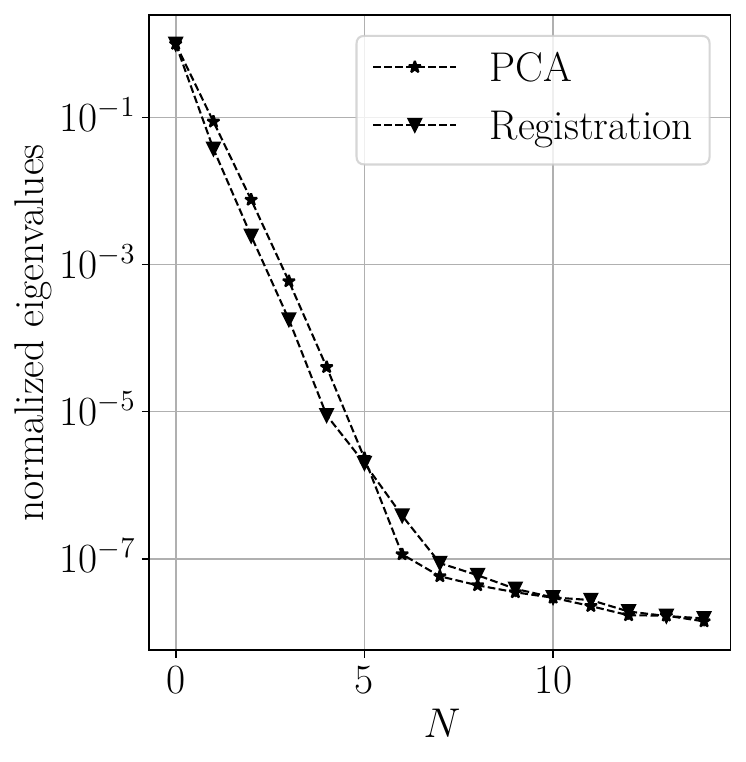}
    }
    \hfill
    \centering\subcaptionbox{\label{subfig:reg_eigenvalues_TD}}{\includegraphics[trim={0cm 0cm 0cm 0cm}, clip, width=0.31\textwidth]{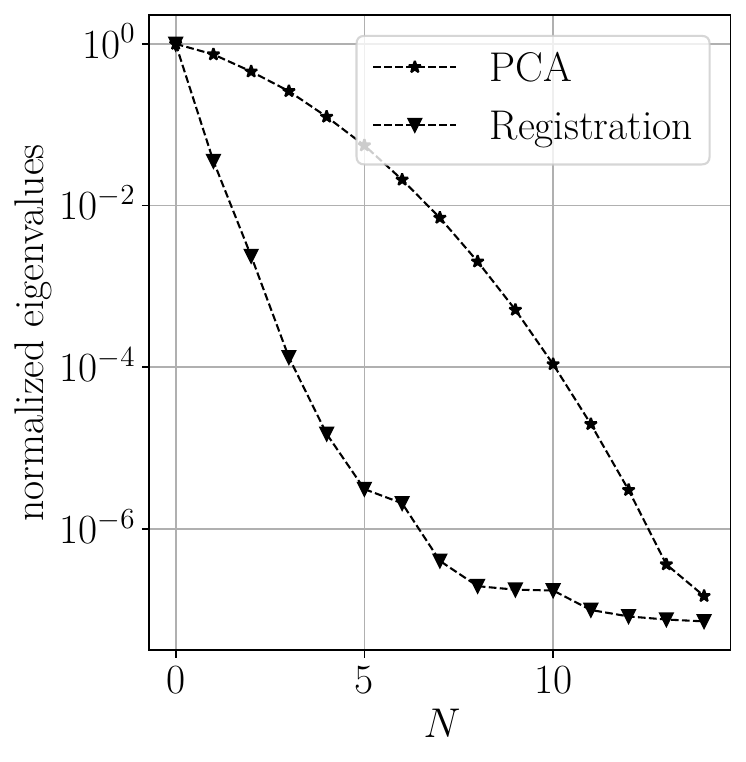}
    }
    \caption{Eigenvalues decay for PCA and Registration methods for pure Advection, pure Diffusion, and Advection-Diffusion problems respectively on the left, in the middle, and on the right.}
    \label{fig:registration_eigen}
\end{figure}

\begin{figure}[b!]
    \centering
    \includegraphics[width=1\linewidth, trim={0cm 4cm 0cm 2cm}, clip]{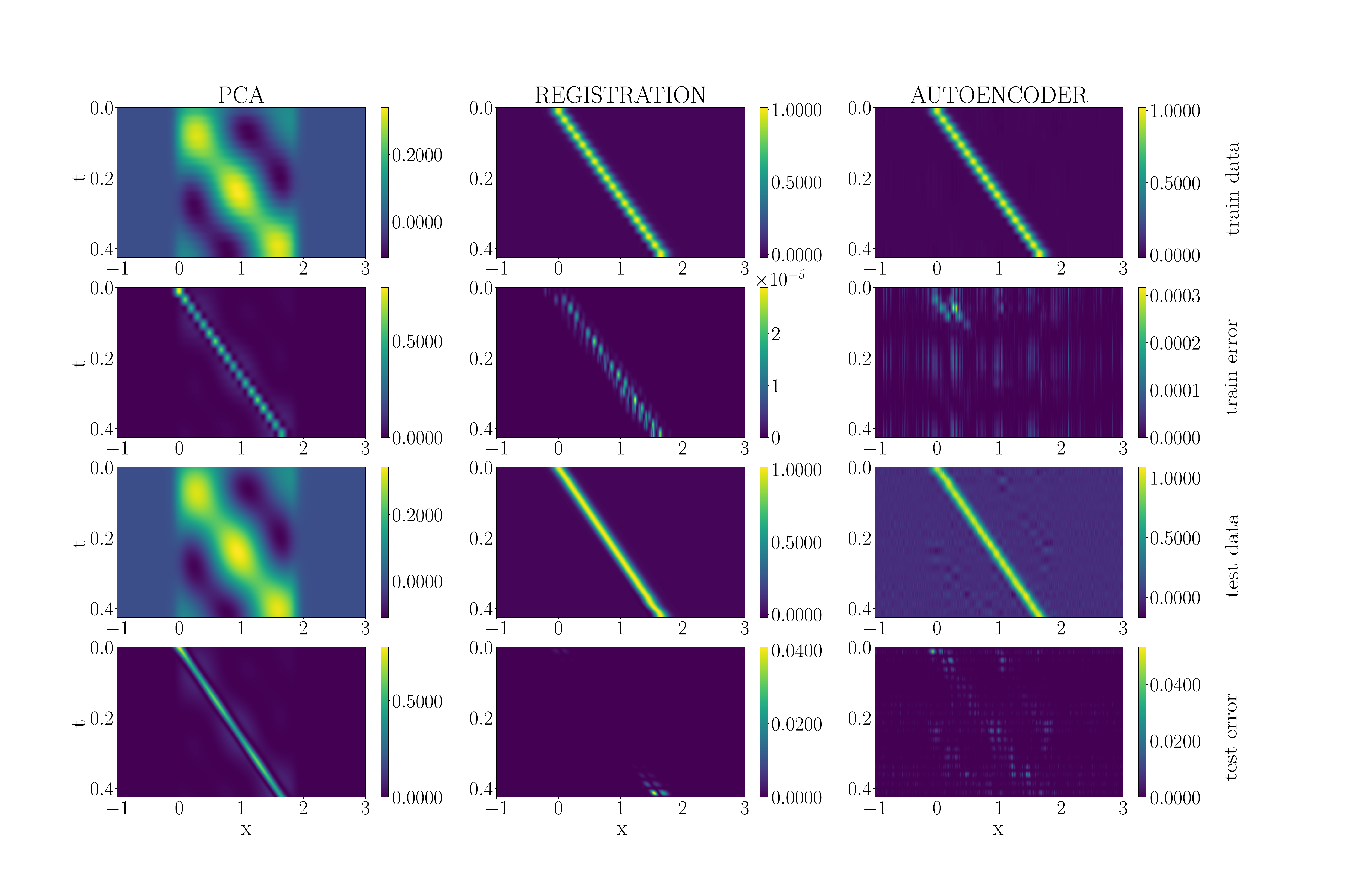}
    \caption{Reconstruction of train and test data shown in first and third row, respectively, with the correspondent $l_2$ errors (second and fourth row) of PCA (left-most column), Registration method (middle column), and AutoEncoder (right-most column) for pure Advection problem with $N=2$.}
    \label{fig:transport_error}
\end{figure}

\begin{figure}[hbt!]
    \centering
    \includegraphics[trim={0cm 4cm 0cm 2cm}, clip, width=1\linewidth]{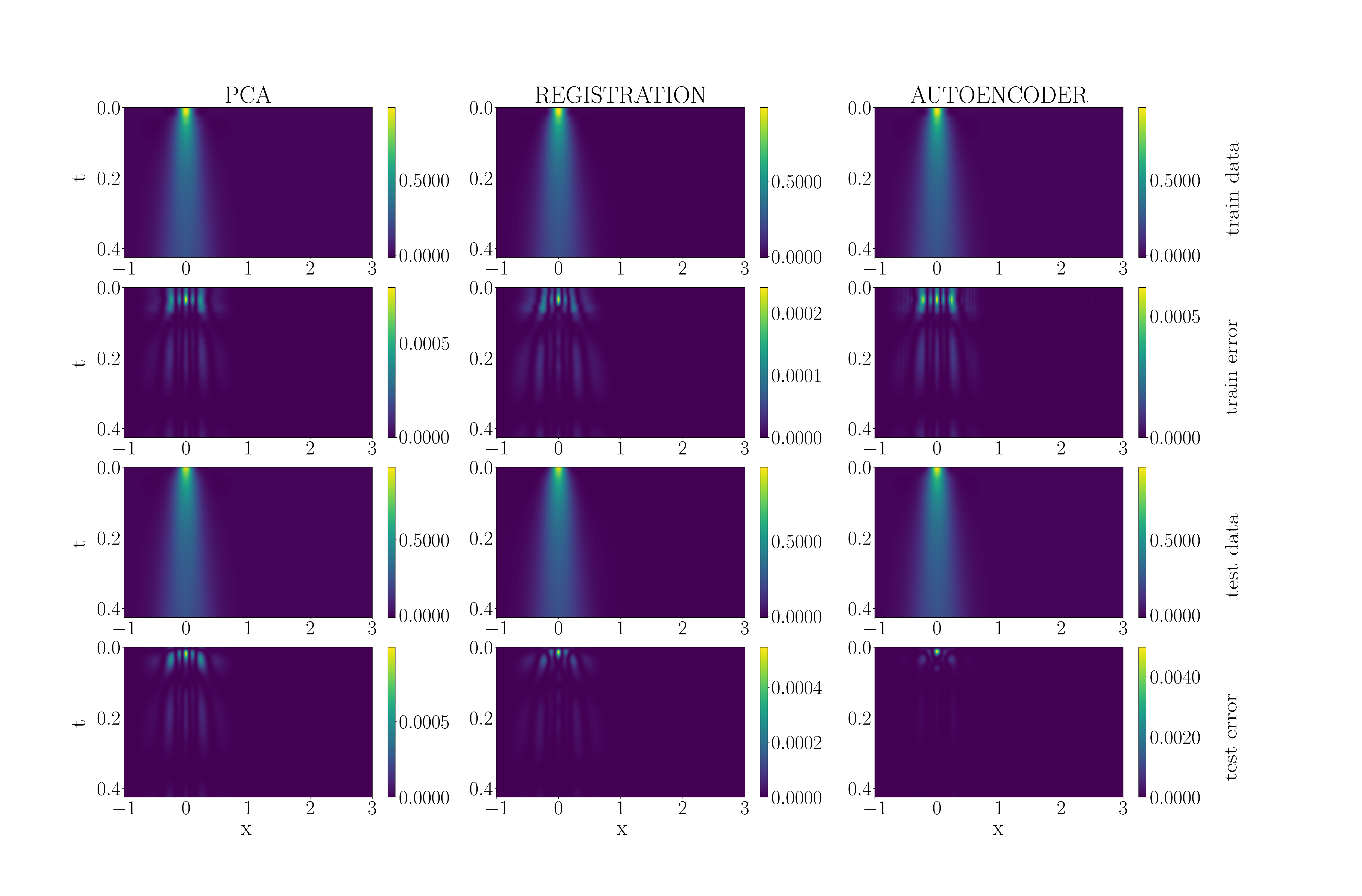}
    \caption{Reconstruction of train and test data shown in first and third row, respectively, with the correspondent $l_2$ errors (second and fourth row) of PCA (left-most column), Registration method (middle column), and AutoEncoder (right-most column) for pure Diffusion problem with $N=2$}
    \label{fig:diffusion_error}
\end{figure}

\begin{figure}[hbt!]
    \centering
    \includegraphics[trim={0cm 4cm 0cm 2cm}, clip,width=1\linewidth]{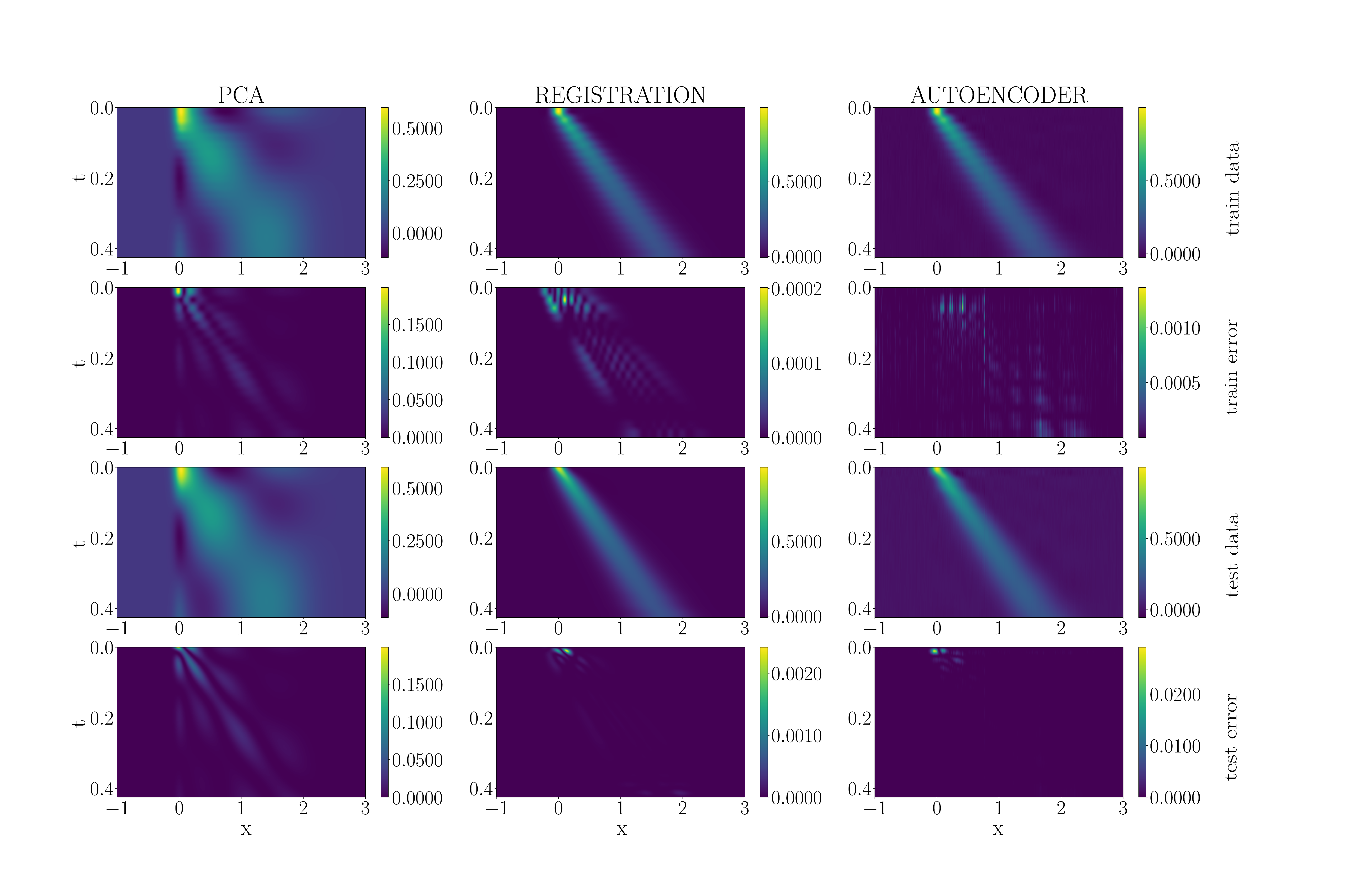}
    \caption{Reconstruction of train and test data shown in first and third row, respectively, with the correspondent $l_2$ errors (second and fourth row) of PCA (left-most column), Registration method (middle column), and AutoEncoder (right-most column) for Advection-Diffusion problem with $N=2$}
    \label{fig:transport_diffusion_error}
\end{figure}

First of all, in Figure \ref{fig:registration_eigen} we report, in the three test cases, the eigenvalues decay of the PCA applied to the dataset $\mathcal{T}$ transformed through the Registration Method, compared with the decay of the PCA performed on the original $\mathcal{M}$ also reported in Figure \ref{subfig:pca_eigen}. As one can see, the Registration procedure appears to be particularly useful in treating Advective behavior, as in both Figures \ref{subfig:reg_eigenvalues_T} and \ref{subfig:reg_eigenvalues_TD} its decay appears much faster than the PCA one.

One can explain this result by recalling the optimization process behind Registration (see Equation \eqref{eq:registration_optimization}). Indeed, it aims to find a suitable map $\Phi$ to bring each datum $u_{t_i}$ as close as possible to a reference one $u_{\overline{t}}$ (here we choose $\overline{t}=0$). This helps in the presence of advective dynamics as it builds a ``less spread" transformed manifold $\mathcal{T}$.

\begin{figure}[b!]
    \subcaptionsetup[figure]{skip=-10pt,slc=off,margin={0pt,0pt}}
    \centering\subcaptionbox{\label{subfig:error_T}}{\includegraphics[trim={0cm 0cm 0cm 0cm}, clip, width=0.31\textwidth]{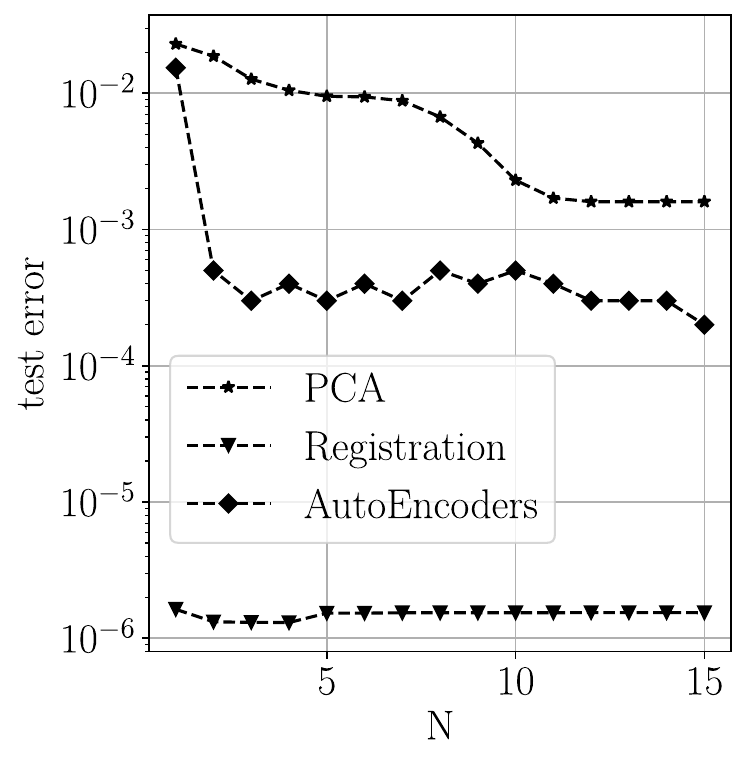}
    }
    \hfill
    \centering\subcaptionbox{\label{subfig:error_D}}{\includegraphics[trim={0cm 0cm 0cm 0cm}, clip, width=0.31\textwidth]{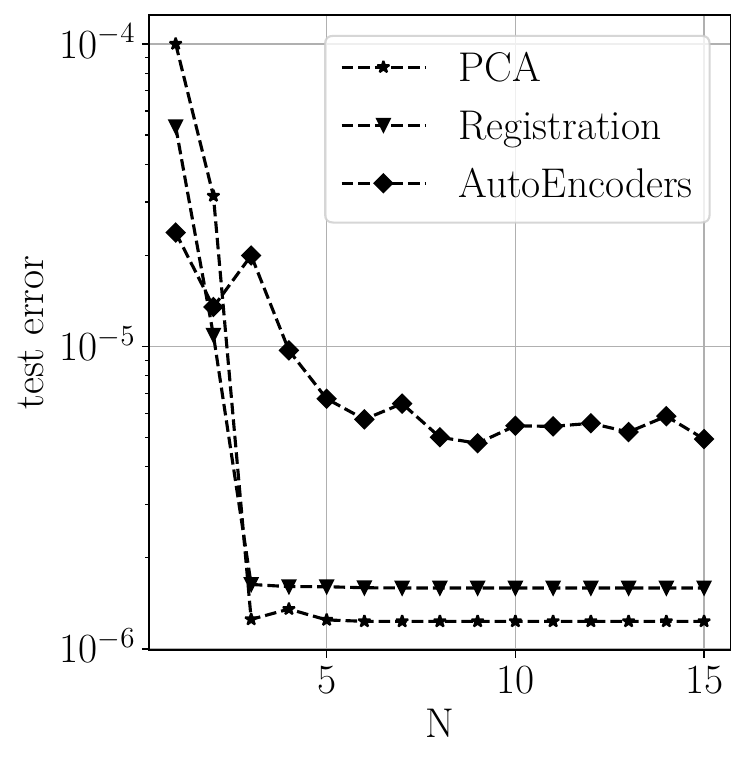}
    }
    \hfill
    \centering\subcaptionbox{\label{subfig:error_TD}}{\includegraphics[trim={0cm 0cm 0cm 0cm}, clip, width=0.31\textwidth]{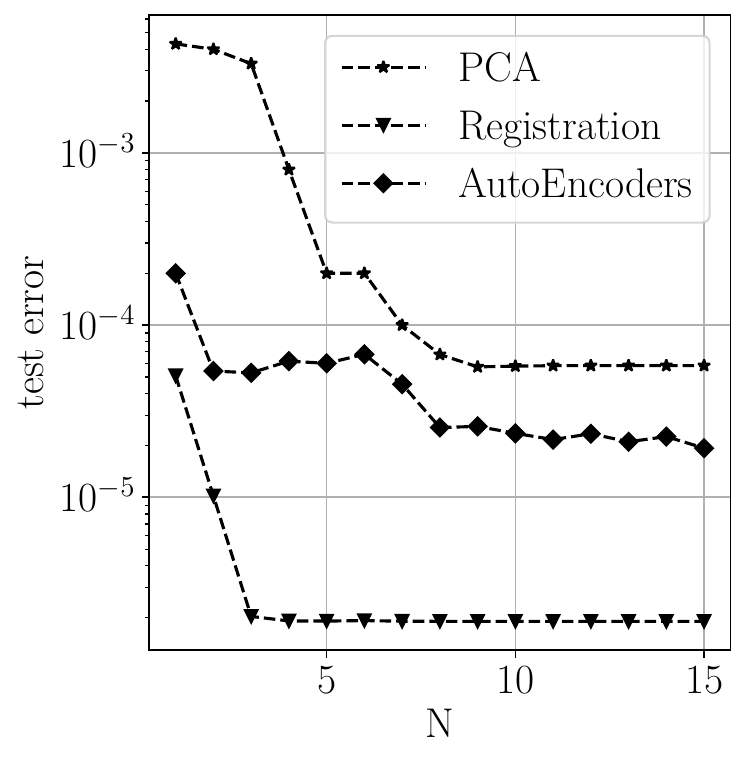}
    }
    \caption{A-posteriori error on the test dataset of PCA, Registration method, and AutoEncoder for pure Advection, pure Diffusion, and Advection-Diffusion problems respectively on the left, in the middle, and on the right.}
    \label{fig:error_reg_AE}
\end{figure}

Furthermore, in Figures \ref{fig:transport_error}, \ref{fig:diffusion_error} and \ref{fig:transport_diffusion_error}, we report the reconstruction with the associated error for the training and testing data, respectively for PCA, Registration methods, and ``vanilla" AutoEncoder for small latent dimension ($N=2$). We construct our ``vanilla" AutoEncoder with Linear Layers interleaved by nonlinear activation function ($\tanh$) for both encoder and decoder, made up also by an hidden layer of dimension $\frac{D}{2}$. 
In all the cases, the test data is reconstructed by applying the backward operators $\Phi^{-1}$ and $Z_N^{-1}$ to the ``test reduced variables'' $z_{\mu_{test}}\in\mathcal{Z}_N$, which are obtained through Kernel Ridge Regression (kRR) based on inverse multiquadratic RBFs \cite{wendland2004scattered} applied on the ``train reduced variables'' $\{\mu_i,z_i\}_{i=1}^n$ collected during the offline phase.
We remark that, in such context, no Manifold Learning strategy is reported due to the fact that in general the map $\Phi$ is not explicitly known, making it impossible to reconstruct the original sample $u_{\mu_i}$ from its latent representation $z_{\mu_i}$. The results highlight the potential of Registration Methods, which, with very few latent variables, provide good approximations in all cases, also comparing to AutoEncoders, while proposing an interpretable reduction pipeline.

\begin{figure}[t!]
    \subcaptionsetup[figure]{skip=-10pt,slc=off,margin={0pt,0pt}}
    \centering\subcaptionbox{\label{subfig:eigenvalues_T}}{\includegraphics[trim={0cm 0cm 0cm 0cm}, clip, width=0.31\textwidth]{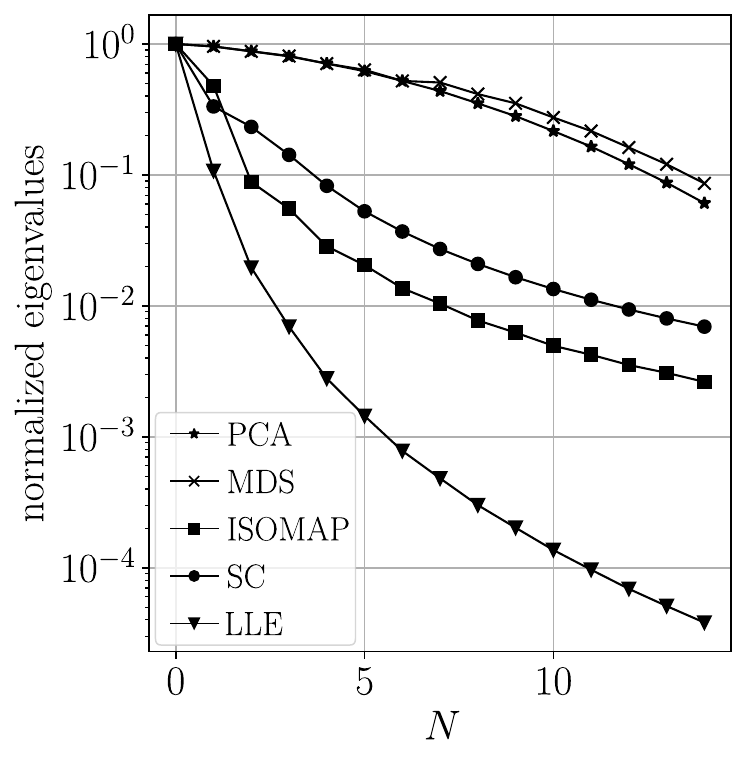}
    }
    \hfill
    \centering\subcaptionbox{\label{subfig:eigenvalues_D}}{\includegraphics[trim={0cm 0cm 0cm 0cm}, clip, width=0.31\textwidth]{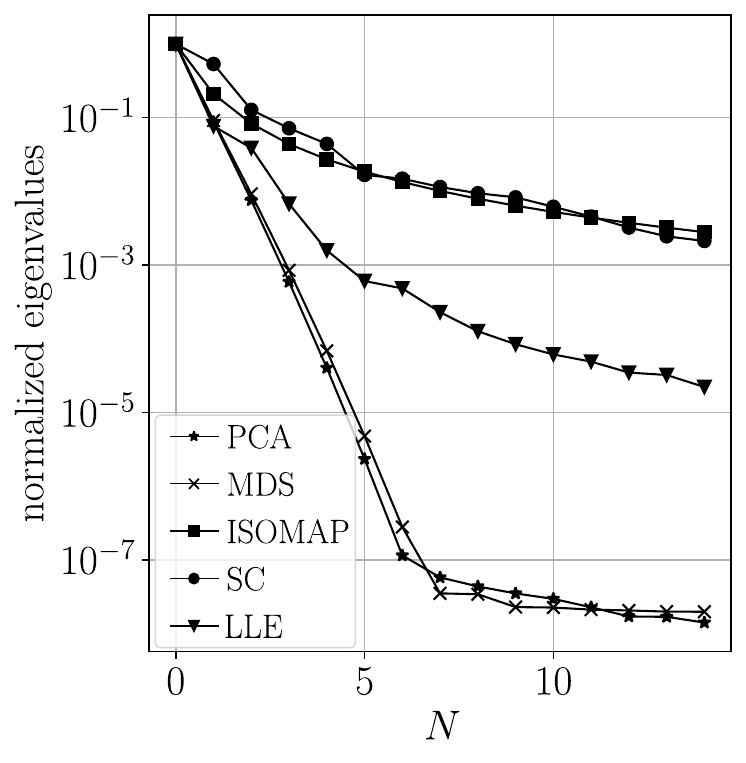}
    }
    \hfill
    \centering\subcaptionbox{\label{subfig:eigenvalues_TD}}{\includegraphics[trim={0cm 0cm 0cm 0cm}, clip, width=0.31\textwidth]{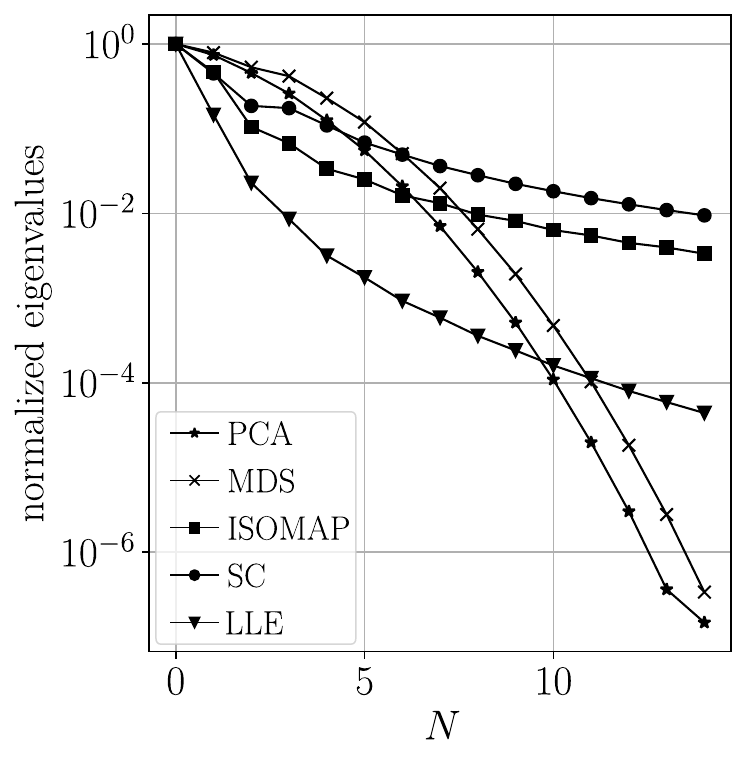}
    }
    \caption{Eigenvalues decay of Manifold Learning methods interpretable as kPCA computed for pure Advection, pure Diffusion, and Advection-Diffusion problems respectively on the left, in the middle, and on the right.}
    \label{fig:KPCA_eigenvalues}
\end{figure}

\begin{figure}[b!]
    \centering
    \includegraphics[width=1\linewidth, trim={0cm 1.6cm 0cm 1.6cm}, clip]{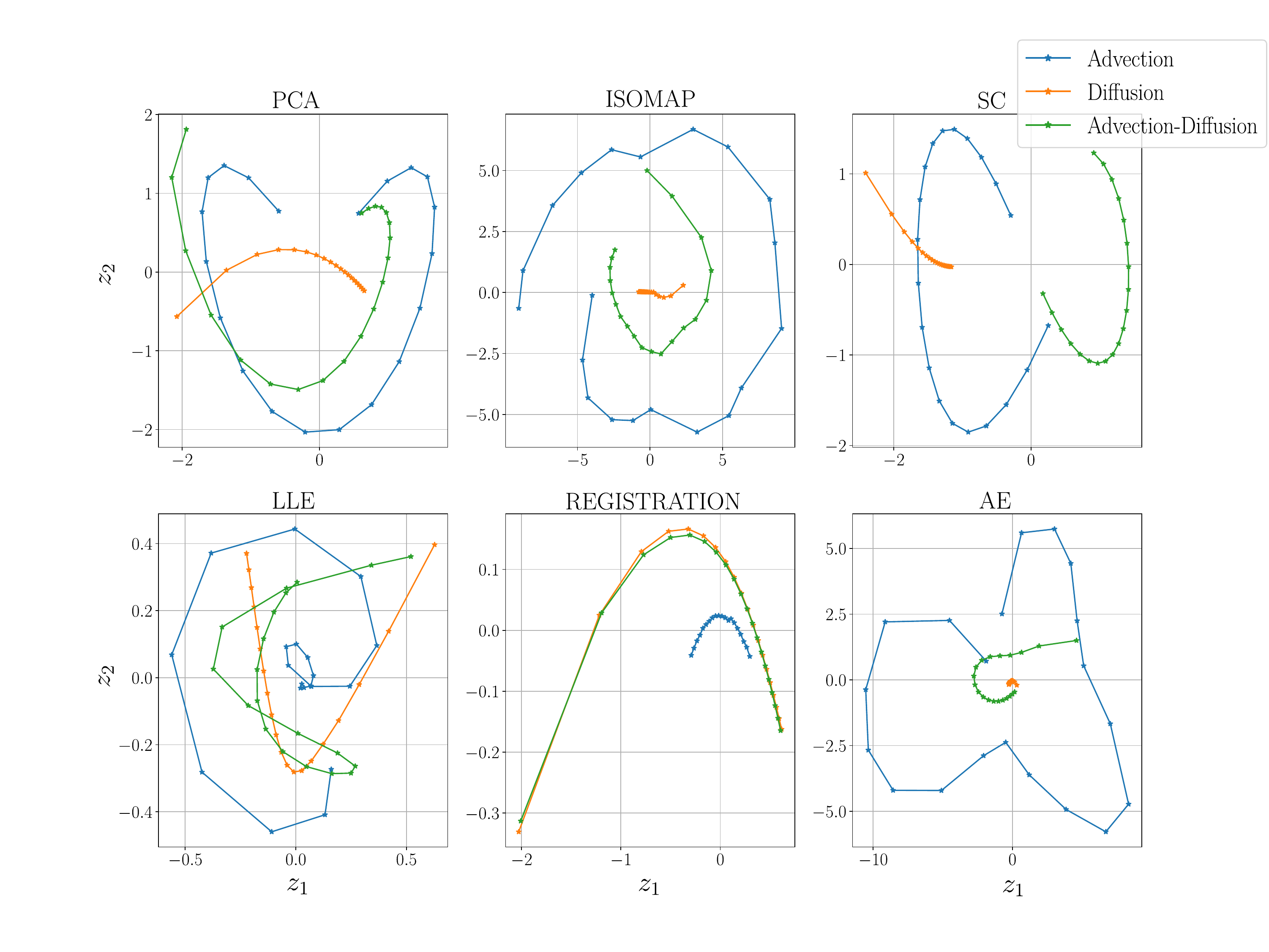}
    \caption{Latent variables dynamic of different nonlinear reduction methods, computed for pure Advection, pure Diffusion, and Advection-Diffusion problems.}
    \label{fig:latent_vars}
\end{figure}

To better capture the error behavior with increasing values of the latent dimension $N$ for both Numerical methods and AutoEncoders, we report in Figure \ref{fig:error_reg_AE} the a-posteriori errors on the test dataset. Regarding PCA and Registration methods, one can note that they decrease as long as their eigenvalues significantly decay, according to Figures \ref{fig:registration_eigen}. Indeed, for Registration the error is stable for $N>5$ in each test case, while the PCA error continues decreasing also for higher $N$ values when the problem is not only Diffusive (Figures \ref{subfig:error_T} and \ref{subfig:error_TD}). Concerning AutoEncoders, their performance does not reach the accuracy of Registration but they demonstrate, in general, an average approximation capability for both Advection and Diffusion problems.

Nevertheless, it is to be underlined that these results correspond to the ``vanilla" forms of the above methods, as many extensions have been proposed in the literature to improve their performance for many test cases \cite{maulik2021reduced,dutta2022reduced,PichiGraphConvolutionalAutoencoder2024,LeeModelReductionDynamical2020,KhamlichOptimalTransportBasedDisplacement2024,KhamlichOptimalTransportinspiredDeep2023}. Since in this work we are interested in analyzing the main ideas behind a reduction process, we do not offer further insights in the direction of ``hyper-optimizing" the architecture, and we only focus on providing a comprehensive analysis of the main strategies giving an intuition on the motivations behind the successful choices for the identification of $\Phi$ and $Z_N$.

Concerning the Manifold Learning techniques, in Figure \ref{fig:KPCA_eigenvalues} we report the eigenvalues decay of the kernel methods of Table \ref{tab:kernel_methods} for the three test cases. For this class of methods, one has to differ between the strategy to adopt for Diffusion and Advection problems in order to obtain the best performance. Indeed, while linear approaches are well-suited when dealing with the former, more advanced techniques are required to obtain a good representation at a low cost for the latter. Above all, LLE seems to provide a good alternative for Advection-dominated problems. This is particularly interesting when considering its similarity to Registration methods. Indeed, both work with mappings that bring one sample into another: it is the case of $\Phi$ such that $u_{\mu_i}\circ\Phi_{\mu_i}\simeq u_{\overline{\mu}}$ for Registration, and of  $\bm{\tilde{W}}$ in Equation \eqref{eq:LLE} for LLE. Therefore, the experiments seem to suggest the utility of considering notions of distance between samples in the analysis of Advection-dominated problems.

Finally, we report in Figure \ref{fig:latent_vars} the reduced variables dynamic for different reduction strategies, i.e.\ their evolution with respect to $t\in[0,T]$ for $N=2$. In general, one can note that the Diffusion problem has a simpler latent dynamic, except that in the subplot corresponding to the Registration method. Indeed, in this case, the dynamics exhibits very similar shapes in all three test cases.
\section{Conclusions}
In this work, we offer a panoramic view of the main ideas concerning nonlinear reduction. In order to do this, we introduce technical insights on various methods belonging to different research fields: Numerical Analysis, Manifold Learning, and Representation Learning.

We highlight their common ground by referring, for their explanation, to the same two-stage pipeline displayed in Figure \ref{fig:general_reduction_scheme}, which we take to be the general scheme of each Nonlinear Reduction strategy. Therefore, interpreting many well-known methods from different fields to the same scheme, we identify their differences by pointing out their diversity in the choices made along the pipeline and explaining their consequences in the effectivity of the final strategy.

Finally, we apply the described methods to the class of Advection-Diffusion problems and analyze their approximation capabilities. In particular, we conduct the analysis on three test cases: pure Diffusion, pure Advection, and Advection-Diffusion problems. We report the reconstruction errors for the classes of Nonlinear ROMs in Numerical Analysis (Registration Methods) and Representation Learning (AutoEncoders), and we analyze the spectral decay for Manifold Learning approaches (kPCA). By analyzing the whole set of approaches, we observed that, as expected, while Diffusion problems seem to be naturally well-suited for linear reduction, Advection problems are successfully reduced only when a proper \textit{pre-processing stage} is applied to the original manifold $\mathcal{M}$. In particular, LLE and Registration Methods seem to lead the way, suggesting to search for a $\Phi$ which transforms all the samples in a reference one or, more in general, that embeds some notion of distance between the training samples.

In conclusion, we foresee that such a common interpretation of both numerical and learning-based reduction can help in future research, by opening broader and deeper communication between the fields.
\section{Acknowledgements}

\textbf{FP} and \textbf{GR} acknowledge the support provided by the European Union- NextGenerationEU, in the framework of the iNEST- Interconnected Nord-Est Innovation Ecosystem (iNEST
 ECS00000043– CUP G93C22000610007) consortium and its CC5 Young Researchers initiative. The
 authors also like to acknowledge INdAM-GNSC and MIUR (Italian ministry for university and research) through FAREX-AROMA-CFD project, P.I. Prof. Gianluigi Rozza, for their support.

\bibliographystyle{abbrv}
\bibliography{bib.bib}
\end{document}